\documentclass[11pt]{article}
\usepackage[noadjust]{cite}
\usepackage{filecontents}
\pdfoutput=1

\usepackage[utf8]{inputenc}
\usepackage{authblk}
\usepackage{indentfirst}
\usepackage{tabularx}
\usepackage{makecell} 
\newcommand{\bhline}{\Xhline{2\arrayrulewidth}}
\usepackage[margin=1.2in]{geometry} 
\usepackage{verbatim}
\usepackage{url}
\usepackage[bookmarks=true, colorlinks=true, linkcolor=blue!50!red, citecolor=blue,
pdfencoding=unicode]{hyperref}
\usepackage[round]{natbib}
\usepackage{subfloat}
\usepackage{caption}
\usepackage{subcaption}
\usepackage{graphicx} 
\graphicspath{{./images/}} 
\usepackage{floatrow} 
\newfloatcommand{capbtabbox}{table}[][\FBwidth]
\usepackage{mathtools} 
\usepackage{amsthm}
\usepackage{amssymb}
\usepackage{fixltx2e} 
\usepackage[T1]{fontenc}
\usepackage{newpxtext}
\usepackage[varg,bigdelims]{newpxmath}
\usepackage[cal=euler,scr=rsfso]{mathalfa}
\usepackage{bm}
\usepackage{inputenc}
\usepackage{microtype}
\usepackage[usenames,dvipsnames]{xcolor}
\usepackage{paralist}
\usepackage{booktabs}
\usepackage{tikz}
\usepackage{makecell}
\usepackage{pgfplots}
\pgfplotsset{compat=1.11} 
\usepackage{todonotes}
\usepackage[capitalize]{cleveref}

\usetikzlibrary{
	cd,
	decorations.markings,
	positioning,
	arrows.meta,
	shapes,
	calc,
	fit,
	quotes}

\usepackage{tabularx} 

\tikzset{
	unoriented WD/.style={
		every to/.style={draw},
		shorten <=-\lsize, shorten >=-\lsize,
		label distance=-2pt,
		thick,
		node distance=\spacing,
		execute at begin picture={\tikzset{
			x=\spacing, y=\spacing}}
		},
	pack size/.store in=\psize,
	pack size = 8pt,
	spacing/.store in=\spacing,
	spacing = \psize,
	link size/.store in=\lsize,
	link size = 2pt,
	pack color/.store in=\pcolor,
	pack color = blue,
	pack inside color/.store in=\picolor,
	pack inside color=\pcolor!20,
	pack outside color/.store in=\pocolor,
	pack outside color=\pcolor!50!black,
	surround sep/.store in=\ssep,
	surround sep=\psize,
	link/.style={
		circle, 
		draw=black, 
		fill=black,
		inner sep=0pt, 
		minimum size=\lsize
	},
	pack/.style={
		circle, 
		draw = \pocolor, 
		fill = \picolor,
		inner sep = .25*\psize,
		minimum size = \psize
	},
	outer pack/.style={
		ellipse, 
		draw,
		inner sep=\ssep,
		color=\pocolor,
	},
	intermediate pack/.style={
		ellipse,
		dashed, 
		draw,
		inner sep=\ssep,
		color=\pocolor,
	},
}

\begin{document}

\title{Evaluating the Pixel Array Method as Applied to Partial Differential Equations}
\author[1]{Cynthia T. Liu}
\author[2]{David I. Spivak}
\affil[1]{\small Massachusetts Institute of Technology, Department of Electrical Engineering and Computer Science}
\affil[2]{Massachusetts Institute of Technology, Department of Mathematics}
\date{}

\maketitle

\begin{abstract}
The Pixel Array (PA) Method, originally introduced by Spivak et. al., is a fast method for solving nonlinear or linear systems. One of its distinguishing features is that it presents all solutions within a bounding box, represented by a plot whose axes are the values of ``exposed variables''. Here we develop a set-theoretic variant of the PA method, named the Pixel Array Solution-Set (PASS) method, that gives PA access to ``hidden variables" whose values are not displayed on plot axes. We evaluate the effectiveness of PASS at numerically finding steady states for several partial differential equations. We discretize several one-dimensional solved reaction-diffusion equations, such as the Fisher equation and the Benjamin-Bona-Mohany equation, using finite differences. Then, we run PASS on each equation, and determine whether it successfully finds all boundary conditions for which a numerical steady state might exist. Then we verify whether the steady states found by PASS are correct. Finally, we discuss the benefits and weaknesses of PASS.
\end{abstract}

\begin{center}
\line(1,0){350}
\end{center}

\section{Introduction}

Many physical systems are modeled using partial differential equations (PDEs), which describe the relationship between changing variables in the system. Examples of PDEs include the reaction diffusion equation, which is used to model chemical reactions, population dynamics, and plasma physics among other physical phenomena; and the Benjamin–Bona–Mahony equation, which models long surface gravity waves \citep{BBM}.

Due to the ubiquity of PDEs, entire fields of research are dedicated to solving them, both analytically and numerically. In particular, numerically solving PDEs has widespread applications in computational fluid dynamics, aerodynamics, and other engineering subjects. 

Equations that have not been solved analytically are often analyzed using numerical methods. For instance, unsolved forms of the reaction diffusion equation:

\begin{equation}
    \label{RD}
    u_t = D(x)u_{xx} + R(u) + f(x,t)
\end{equation}

Are examined using numerical methods in \citep{Numerical_RD_1, Numerical_RD_2, Numerical_RD_3}. These methods usually require initial and boundary conditions to be satisfied, and provide one numerical solution.

This is where the Pixel Array (PA) method \citep{Introduction_to_PA} becomes relevant. Unlike other numerical methods, it provides all solutions to a system of equations within a bounding box and is often faster than quasi-Newton methods. Because partial differential equations can be converted to a system of equations using discrete approximation methods, the PA method can be applied to partial differential equations. However, the PA method is quite new, and thus the goal of this paper is to establish its accuracy at numerically solving PDEs.

To test and also use the PA method, we needed to adapt it so as to produce a matrix of full solution sets, rather than just boolean values indicating where solutions exist. We refer to this adapted method as the Pixel Array Solution-Set (PASS) method. The ``bounding box" in this case consists of a pair of ranges for boundary values, one on each side of the one-dimensional PDE system; of course a higher-dimensional tensor of solution sets would be appropriate for a higher-dimensional system. Using the PASS method, we were able to accurately find steady states for a variety of well-established PDEs.

The plan of the paper is as follows. In \cref{sec:methods} we provide a brief review of the finite differences method, as well as describe our modification of the PA method. In \cref{sec:testing} we use finite differences to discretize the heat equation, the local Fisher equation, the Benjamin-Bona-Mahony equation, and the Sine-Gordon equation. Then we use PASS to calculate numerical solutions. As all the equations were previously solved, we can use them to demonstrate the effectiveness of PASS.

\section{Methods}\label{sec:methods}

The methods we use to numerically solve for steady states are described below. We first briefly illustrate our discretization method, finite differences. Then we describe the difference between the Pixel Array method of \citep{Introduction_to_PA} and the Pixel Array Solution-Set (PASS) method developed here. Finally, we demonstrate how using PASS with discretized partial differential equations lets us generate steady states.

\subsection{Finite Differences}

To apply the PASS method to a PDE, we need to remove all derivatives from the equation. We choose to do this with finite differences, which uses secant lines to approximate derivatives. Using finite differences, the approximate second derivative with respect to $x$ is:

\begin{equation}
    \label{FD second order}
    u_{xx} = \cfrac{1}{h^2} (u_{i+1} - 2u_i + u_{i-1})
\end{equation}

\noindent Where $h$ is the length of a discrete spatial interval in $x$, and $u_i$ corresponds to numerical solution value in the spatial interval $i$. As an example, we apply finite differences to the reaction-diffusion equation \eqref{RD}. Note that because we only care about steady states, any derivatives with respect to time $t$ are set to 0. We get:

\begin{equation}
    \label{discrete RD}
    0 = \cfrac{D}{h^2} (u_{i+1} - 2u_i + u_{i-1}) + R(u_i) + \hat{f}(x_i)
\end{equation}

Equations like \eqref{discrete RD} will be called \textit{discrete steady state conditions}. We will be solving systems of that equation for all $i$ with the PASS method. If the discrete steady state condition is satisfied for all $u_i$, then the values of $u_i$ $\forall i$ in order form a numerical steady state.

\subsection{PASS: Modifying the Pixel Array method to pass along solution sets}\label{sec:PASS}

\begin{figure}[h]
\begin{subfigure}{.4\textwidth}
  \centering
  \captionsetup{width=0.8\textwidth}
  \includegraphics[width=\linewidth]{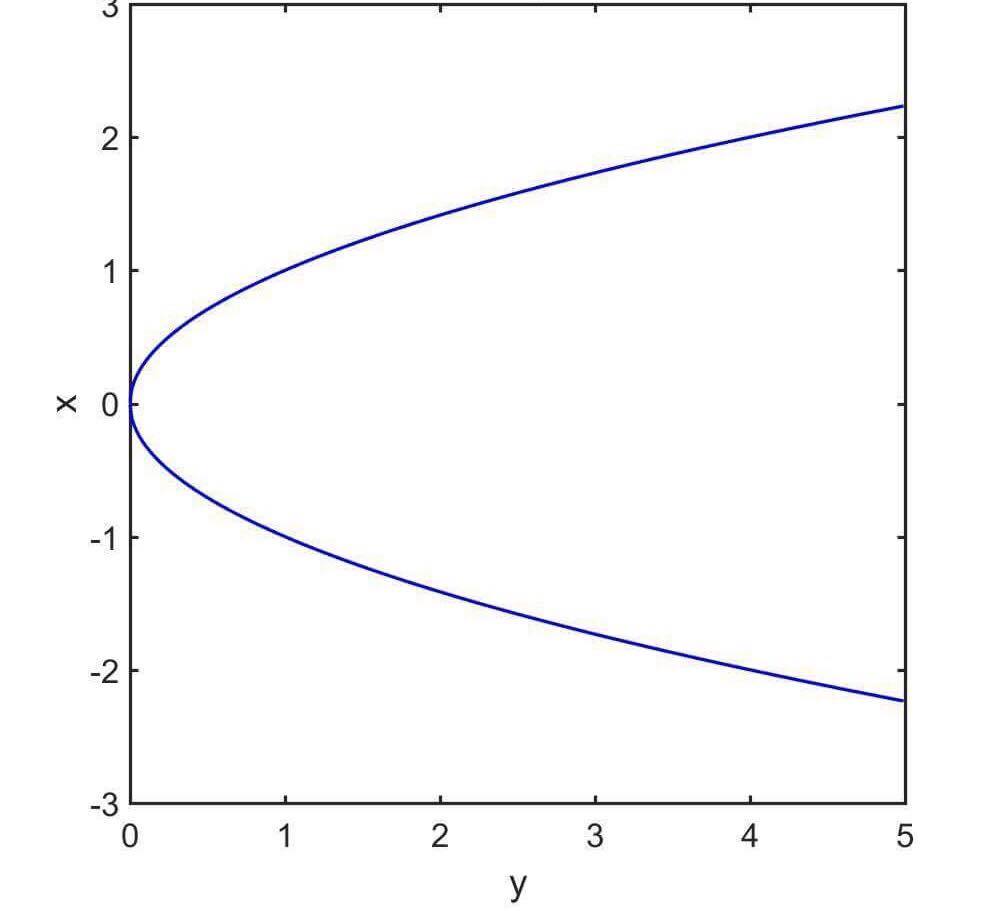}
  \caption{An example plot of a function, oriented so its axes match that of a matrix}
  \label{sample_function}
\end{subfigure}%
\begin{subfigure}{.4\textwidth}
  \centering
  \captionsetup{width=0.84\textwidth}
  \includegraphics[width=.84\linewidth]{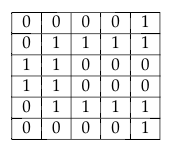}
  \caption{A low-resolution pixel array corresponding to the function plot.}
  \label{sample_pixel_array}
\end{subfigure}%
\caption{An example plot with its corresponding pixel array}
\label{plot_and_pa}
\end{figure}

Each equation in the discretized system corresponds to an array in both the Pixel Array (PA) and PASS method. In \citep{Introduction_to_PA}, the PA method is executed with boolean k-dimensional matrices or tensors, each of which represents one equation or inequality. Each of the $k$ dimensions in a given tensor specifies one variable of the corresponding equation, and each row or column in that dimension represents an interval of values the variable can assume. For each possible combination of $k$ variable values, the corresponding entry in the matrix is ``True" (1) if the equation is satisfied somewhere in that range of values and ``False" (0) otherwise. To represent this, we present Figure \cref{sample_function}, which is a function graph, and figure (\cref{sample_pixel_array}) as its corresponding pixel array with vertical dimension y and horizontal dimension x. Observe that the plot and pixel array are turned to their sides, to match conventions described in \citep{Introduction_to_PA}.

In Figure (\ref{plot_and_pa}), each of the cells in Figure (\ref{sample_function}) corresponds to a pixel in Figure (\ref{sample_pixel_array}). For instance, the function $y=x^2$ goes through the cell with $x \in [1,2]$, $y \in [2,3]$, and therefore its associated pixel, the third column in the second row, has the value 1.

It is useful to set up a naming convention which relates the tensor entries and the parts of the plot they correspond to. We thus refer to a tensor entry by the value exactly at the center of its range; e.g. in the above example our pixel would be $(1.5, 2.5)$. When we wish to refer to all of an entry's possible values as opposed to just the center value, we refer to it as a \textit{subcube}. We say that the \textit{bin-size}---the range of each variable defining subcube---for this pixel is $(1,1)$. For simplicity, we will assume that the bin sizes of all variables are equal, for every system we discuss. 

Suppose each equation in the system has been plotted as a tensor, as above. To solve the simultaneous system of equations, these tensors are multiplied together by contracting various edges. The resulting tensor will then indicate all combinations of the remaining variables for which there exists a solution to the system.

What differentiates the PASS method from the PA method is that we are interested in finding numerical steady states, not just determining their existence for certain boundary conditions. Thus, the boolean definition of a pixel array is too limited. This is resolved by introducing a new data type, the \textit{solutionSet}. Each entry of our modified pixel array is of this data type. A solutionSet is a finite set of tuples of any finite length. In our work, the lengths of each tuple in a given solutionSet turn out to be the same, e.g.\  $\{(0,1,2), (1,1,1), (2,1,0)\}$, because they will represent the set of solutions---in certain selected variables---to some system of equations.

To get some intuition for the solutionSet data type, consider the homogeneous heat equation:

\begin{equation*}
    \label{general homogeneous heat}
    u_t - D(x)u_{xx} = 0
\end{equation*}

\noindent With $u_t = 0$ and $D(x) \equiv 1$, the equation becomes

\begin{equation}
    \label{homogeneous heat}
    u_{xx} = 0,
\end{equation}

\noindent otherwise known as the one-dimensional Laplace equation. After being processed by finite differences and multiplying both sides by $h^2 \ne 0$ we get: 

\begin{equation}
    \label{discrete heat}
    u_{i+1} - 2u_i + u_{i-1} = 0
\end{equation}

\noindent The corresponding tensor has three dimensions, corresponding to the variables $u_{i+1}$, $u_i$, and $u_{i-1}$ in \cref{discrete heat}. Each tensor entry corresponds to a subcube; let $(x,y,z)$ denote the center point of that subcube. The value in that entry is the solution-set for the equation \cref{discrete heat}. That is, the value will be the empty set $\{\}$ if the subcube centered at $(x,y,z)$ does not contain a solution $(u_{i+1}, u_i, u_{i-1})$ where value of each $u_j$ is within the subcube, and it will be the one-element set $\{(x,y,z)$\} if there exists such a $(u_{i+1}, u_i, u_{i-1})$.

Tensor multiplication involves multiplication and addition. For boolean values, multiplication is as usual, e.g.\ $1*0=0$, and addition is almost as usual, the exception being that $1+1=1$; this is the usual \emph{Boolean semiring}, which is used in the PA method. For the PASS method, we need to define multiplication and addition of the solutionSet data type. This data type will store the solutions to the local discrete steady state condition.

The multiplication of solutionSets is defined to be the Cartesian product of the tuples, and addition is defined to be set union. So if two theoretical solutionSets, \{(1),(2),(3)\} and \{(3,4)\} were multiplied, the resultant tuple would be \{(1,3),(2,3),(3,3),(1,4),(2,4),(3,4)\}. If they were added, we would get \{(1),(2),(3),(4)\}. This defines a commutative semiring, whose additive unit is the empty set $\{\}$ and whose multiplicative unit is $\{()\}$. Thus for example, when a solutionSet is multiplied by $\{\}$, the result is zero. When a solutionSet is added to $\{\}$ or multiplied by $\{()\}$, the solutionSet remains the same.

With this definition of multiplication, the generalized matrix multiplication algorithm remains the same as it was in \citep{Introduction_to_PA}.

\subsection{Producing Solutions}

We will now demonstrate how the above definition of the solutionSet data type allows us to compute numerical steady states to a partial differential equation. Specifically, we will show how the numerical steady states are stored within solutionSets, and expanded as more multiplications are performed. 

We show this by example. Once again, consider the heat equation \eqref{homogeneous heat} and its discrete version \eqref{discrete heat}. As discussed in \cref{sec:PASS}, all entries in the pixel array with dimensions $u_{i+1}$, $u_i$, and $u_{i-1}$ contain the single solution ${(u_i)}$ if $u_i$ satisfies \cref{discrete heat}. This is true for all partitions $i$.

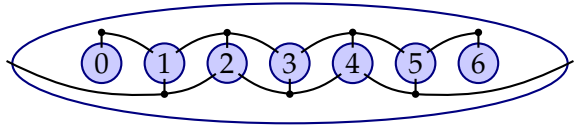
\begin{figure*}
\caption{Wiring diagram for the discrete heat equation with n = 7 partitions. The boundary conditions in this instance are cells 1 and 5, but the product is evaluated as if 0 and 6 were boundaries to help prevent discontinuities at the true boundary.}
\label{WD_heat}
\begin{center}
\begin{tikzpicture}[unoriented WD, surround sep= 3pt, every to/.style={draw, bend right=15pt}]
	\node[pack]                 (n0) {0};
	\node[pack, right=1 of n0]  (n1) {1};
	\node[pack, right=1 of n1]  (n2) {2};
	\node[pack, right=1 of n2]  (n3) {3};
	\node[pack, right=1 of n3]  (n4) {4};
	\node[pack, right=1 of n4]  (n5) {5};
	\node[pack, right=1 of n5]  (n6) {6};
	\node[link, above=.3 of n0] (l0) {};
	\node[link, above=.3 of n2] (l2) {};
	\node[link, above=.3 of n4] (l4) {};
	\node[link, above=.3 of n6] (l6) {};
	\node[link, below=.3 of n1] (l1) {};
	\node[link, below=.3 of n3] (l3) {};
	\node[link, below=.3 of n5] (l5) {};
	\node[outer pack, fit=(n0.center) (l2) (l3) (n6.center)] (outer) {};
	\draw (n0.north) -- (l0);
	\draw (n2.north) -- (l2);
	\draw (n4.north) -- (l4);
	\draw (n6.north) -- (l6);
	\draw (n1.south) -- (l1);
	\draw (n3.south) -- (l3);
	\draw (n5.south) -- (l5);
	\draw (n1) to (l0);
	\draw (n2) to (l3);
	\draw (n3) to (l2);
	\draw (n4) to (l5);
	\draw (n5) to (l4);
	\draw (l1) to (n2);
	\draw (l2) to (n1);
	\draw (l3) to (n4);
	\draw (l4) to (n3);
	\draw (l6) to (n5);
	\draw (outer.west) to (l1);
	\draw (l5) to (outer.east);
\end{tikzpicture}
\end{center}
\end{figure*}

With all the matrices defined, we will multiply them using the general matrix multiplication algorithm described in \cref{sec:PASS}, but to do so we must determine which variables are exposed. We resolve this by introducing the \textit{Wiring Diagram} (WD) for this system; see \cref{WD_heat}. As described in \citep{Introduction_to_PA}, the WD lets us connect matrices by their variables and lets us decide which variables we want to be exposed when we have finished multiplying.

In \cref{WD_heat}, the black dots (and the lines connecting them to circles) are called \textit{links}, while the circles are called \textit{cells} (in \citep{Introduction_to_PA} they were called \textit{packs}). Each cell represent one equation, its ports represent the variables of that equation, and the links represent the dependencies between variables. For example, because the equation for $u_i$ in \cref{discrete heat} depends on $u_{i+1}$ and $u_{i-1}$, there are links connecting the cells representing $u_{i+1}$, $u_{i-1}$, and $u_i$.

The last thing to explain about \cref{WD_heat} is the fact that we have both \emph{hidden boundaries}, cells 0 and 6, and the \emph{exposed boundaries}, cells 1 and 5. The boundary conditions for the heat equation reference the exposed boundaries. At first, we did not include hidden boundaries in our research. However, we found that in that case, we get many discontinuous or anomalous steady states due to the spacial discretization of equations with a left-right asymmetry. For example, this problem does not arise for the heat equation, because our equation is symmetric (it treats $u_{i+1}$ and $u_{i-1}$ the same way). By adding the hidden boundaries, the anomaly is partially handled.

\begin{figure*}
\caption{Wiring diagram for a single product with 4 variables, all exposed}
\label{WD_oneproduct}
\begin{center}
\begin{tikzpicture}[unoriented WD, surround sep= 3pt]
	\node[pack] (n1) {$M_1$};
	\node[pack, right=1 of n1] (n2) {$M_2$};
	\node[link, above=.3 of n2] (l2) {};
	\node[link, below=.3 of n1] (l1) {};
	\node[outer pack, fit={(n1) (n2) (l1) (l2)}] (outer) {};
	\draw (n1.south) -- (l1);
	\draw (n2.north) -- (l2);
	\draw (n1.north west) to[bend left=15pt] (outer.north west);
	\draw (n2.south east) to[bend left=15pt] (outer.south east);
	\draw (l1) to[bend right] (n2.south west);
	\draw (l2) to[bend right=15pt] (n1);
	\draw (l1) to[bend left=15pt] (outer.west);
	\draw (l2) to[bend left=15pt] (outer.east);
\end{tikzpicture}
\end{center}
\end{figure*}
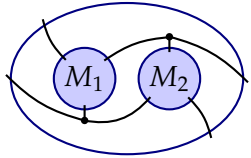

\cref{WD_oneproduct} represents a multiplication of two 3-dimensional tensors $M_1$ and $M_2$, the result of which is a 4-dimensional tensor.\footnote{The dimension of a tensor is indicated in the diagram by its number of ports.} It is given by contracting along the links connecting $M_1$ and $M_2$, which represent shared dimensions. Performing multiple such tensor multiplications, one obtains the tensor multiplication indicated by \cref{WD_heat}.


Now that we have explained how to multiply matrices and determine exposed variables and boundary conditions, we proceed to compute the value of a general solutionSet.  By the definition of solutionSet multiplication defined in Section 2.2, if the entry $\{(u_1)\}$ is multiplied by the entry $\{(u_2)\}$, the resulting solutionSet is $\{(u_1, u_2)\}$. Recall from \cref{sec:PASS} that the entry $\{(u_1)\}$ would only occur if \cref{discrete heat} were satisfied for that particular value. The same can be said for the entry $\{u_2\}$. Therefore, the product $\{(u_1, u_2)\}$ would only be produced if their values satisfied in \cref{discrete heat} for \textit{both} $i = 1$ and $i = 2$. This means that the final value within the solutionSet is a tuple of numerical steady states, because each local numerical solution fulfills the discrete steady state condition.

Thus, by multiplying all the matrices together, we calculate solutionSets of numerical steady states. Each element within a single solutionSet satisfies the discrete steady state condition locally. Furthermore, with the PASS method, we can calculate steady states for all possible pairs of left and right boundary values at once, in some sense parallelizing that computation.

\section{Testing}\label{sec:testing}

We analyze the effectiveness of PASS by using it to produce steady states for 4 different partial differential equations. In \cref{sec:homog_heat} we analyze the heat equation. In \cref{sec:Fisher} we test on the Fisher equation. Then, in \cref{sec:bbm}, we examine the Benjamin-Bona-Mohany equation. Finally, in \cref{sec:sg}, we look at the Sine-Gordon equation.

\subsection{Homogeneous Heat Equation}\label{sec:homog_heat}

We first consider the homogeneous heat equation as described by \cref{homogeneous heat}, whose discrete steady state condition is \cref{discrete heat}. It is well known that steady states to the homogeneous heat equation are linear between the two boundary conditions, and this fact can also be verified by directly solving the equation. 

Let us first check if PASS successfully finds all boundary conditions for which there are steady states. For the heat equation, any pair of boundary conditions should have an associated steady state. For possible boundary values in the range $[0,2]$, our associated pixel array is:

\begin{figure}[h]
\centering
\includegraphics[width=4.5cm]{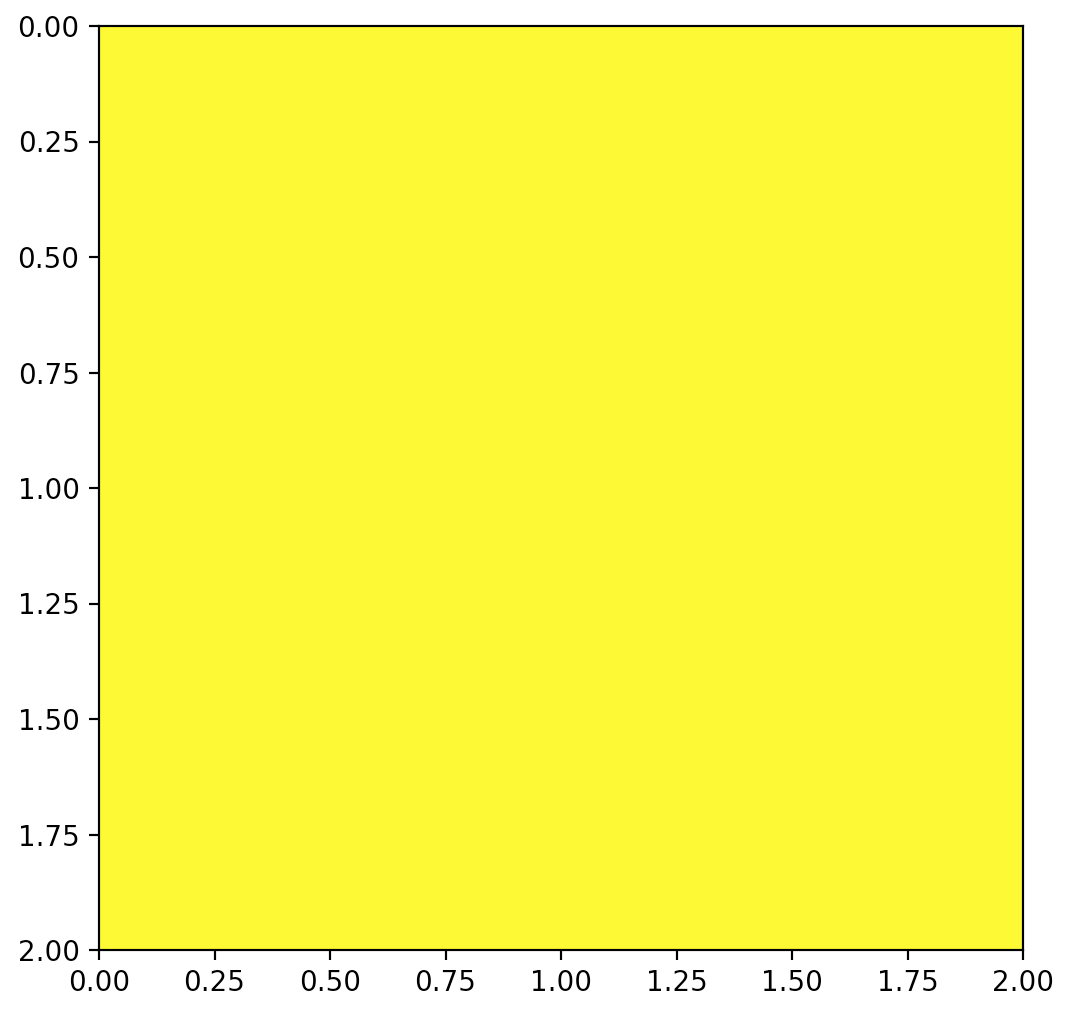}
\caption{The pixel array associated with the heat equation. Yellow indicates that at least one numerical steady state was found for the boundary conditions indicated by the axes.}
\label{heatPixel}
\end{figure}

This is our desired result. Now let us check whether the numerical steady states found by PASS are sufficiently similar to the true, perfectly linear, steady state. One can prove that, for any number of cells and any range, as the number of bins tends to infinity, the numerical solution given by the PASS method converges to the true solution. Indeed, this is a direct consequence of Theorem 3.4.1 in \citep{Introduction_to_PA}.
%
%
%
%

With a finite number of bins, we can also theoretically bound the 'distance' between the true solution and the results of the PASS method. Suppose there are $n$ cells, and for each cell $1\leq i\leq n$, let $u_i$ be the true solution to the heat equation at $x_i$ (these should be a linear function of $i$) and let $v_i$ denote the numerical solution. By the definition of bins and the $L^\infty$ norm, the coordinates at any vertex of a subcube are $\pm \frac{b}{2}$ away from the center coordinates. 

For a numerical solution to the heat equation to exist within the subcube, the value of equation \cref{discrete heat} must be within $\pm \epsilon$ of 0 for some $\epsilon > 0$. We calculate epsilon by finding the maximum magnitude of the gradient within the subcube, and then multiplying that by the ``radius,'' or largest distance to any boundary, of the subcube. This calculation gives $\epsilon = \frac{3b\sqrt{2}}{2}$. Thus the maximum difference from 0 at the center that still permits some point in the subcube to have the value 0 is $\frac{3b\sqrt{2}}{2}$. In general, plotting a pixel array by comparing an expression's value at the center of a subcube to an $\epsilon$ will be called the \emph{epsilon method}.

We can also calculate a functional bound of how much our solution could differ from the true steady state. We define this as the square root of the sum of the squared differences $(v_i-u_i)^2$ for all i. After some calculation, we get:

\begin{equation}
\label{n norm}
L = \epsilon \left \lceil \cfrac{n}{2} \right \rceil^{5/2}
\end{equation}
where $L$ is the \emph{modified $L_2$ norm}. $L$ clearly goes to 0 as $b$ goes to 0. 

\begin{figure*}
\begin{center}
\begin{tikzpicture}
\begin{axis}[
title = \footnotesize Numerical Approximations to Steady States,
xlabel = \footnotesize Cell Number $i$,
ylabel = \footnotesize Value at Cell,
xmin = 1, xmax = 8,
width=2.5in
]
\addplot[black, mark=o] coordinates {
	(1, 1.0) 
	(2, 1.0)
	(3, .95)
	(4, .85)
	(5, .7)
	(6, .55)
	(7, .35)
	(8, .1)
};
\addplot[blue, mark=square] coordinates {
    (1, 1.0)
    (2, .85)
    (3, .62)
    (4, .45)
    (5, .3)
    (6, .2)
    (7, .15)
    (8, .1)
};
\addplot[red, mark=triangle] coordinates {
    (1, 1.0)
    (2, .95)
    (3, .85)
    (4, .7)
    (5, .55)
    (6, .45)
    (7, .3)
    (8, .1)
};
\end{axis}
\end{tikzpicture}
\caption{Some numerical steady state approximations given for homogeneous heat equation with the boundaries [1.0,.1]}
\label{heats}
\end{center}
\end{figure*}
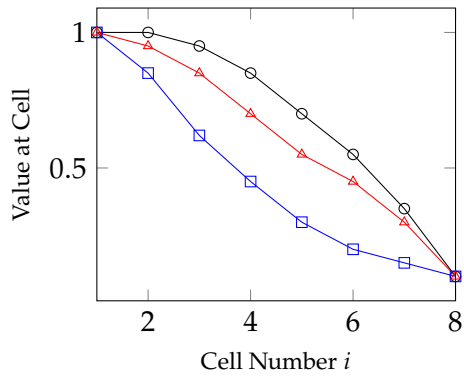

We run the PASS method%
\footnote{The reader can find the code for this method at \url{https://github.com/cynliu98/Pixel\_Array\_Python}, called \textit{PASS.py}, in the archives folder.}
for a few combinations of $(b,n)$ to verify the bound in \eqref{n norm}. The program partitions the spatial dimension into $n=8$ parts (which we will call partitions) and a temperature mesh of $b = .05$, giving 21 bins for each cell. Using the single boundary value pair $[1.0,0.1]$, the PASS method generated a total of 10 steady state approximations. 3 of those numerical steady states generated by PASS are shown in \cref{heats}. As one can see, they are all approximately linear---a perfectly linear steady state is the true steady state of the the equation.

The largest modified $L_2$ norm among PASS's solutions was 1.46. When $b = .05$ and $n = 8$ are substituted in equation $\cref{n norm}$, we get $3.39$, which is greater than $1.46$. Thus, all steady states' norms are smaller than the bound determined by equation \cref{n norm}. This bound also holds for the steady states given by other pairs of boundaries. Furthermore, in figure \cref{heats}, we see that the numerical steady states are approximately linear, giving visual confirmation that our approximations are close to the steady states.. The runtime under these conditions was about 48 seconds in Python 3 for 441 boundary value pairs on a 2016 Macbook. 

\subsection{Fisher-KPP Equation}\label{sec:Fisher}

Now, we consider the classical Fisher-KPP equation \citep{Fisher}. In one dimension, the equation is:

\begin{equation}
    \label{Fisher}
    u_t = u_{xx} + \mu u (1-u)
\end{equation}

\noindent Where $\mu$ represents the intrinsic growth rate of a species \citep{Numerical_RD_1} when used in the field of population dynamics. Its corresponding discrete equation for finding steady states is:

\begin{equation}
    \label{discrete Fisher}
    0 = \cfrac{1}{h^2} (u_{i+1} - 2u_i + u_{i-1}) + \mu u_i (1-u_i)
\end{equation}

\noindent The only bounded steady states of this equation are $u \equiv 1$ and $u \equiv 0$ \citep{Fisher}, for all values of $\mu > 0$ and spatial interval size $h > 0$.

Our method of plotting the pixel arrays for this equation is slightly different than plotting the array for the heat equation. It is essentially a stricter form of the epsilon method as defined in the previous section. An outline of the method is:

\begin{enumerate}
    \item Let the discrete reaction diffusion expression be $f(u_{i+1}, u_i, u_{i-1})$, and calculate its value.
    \item  Round that expression to the nearest bin, and round 0 to the nearest \emph{expanded bin}, where possible bin values $b_j$ are written as $b_j = mj + c$ for $j \in \mathbb{Z}$, $j \in [0,r]$ where $r$ is the resolution, and expanded bin values are just values of $b_j$ with unbounded, integer $j$. If those rounded values are the same, then the corresponding entry in the Pixel Array is nonempty.
\end{enumerate}

This method obviously converges as the size of bins $b \rightarrow 0$, because $f(u_{i+1}, u_i, u_{i-1})$ and 0 must round to the same expanded bin and smaller $b$ restricts the range of values that could round to the same bin. This is true for any fixed spatial interval size $h$, including as $h \rightarrow 0$.

We run the program with 41 bins, $b = .05$ in the range $[0,2]$ and 10 cells, rounding up. We first test 3 values of $h$, $h = 1$, $h = .25$, and $h = .1$, all with $\mu = 1$. Because there are far significantly fewer steady states (just $u \equiv 1$ and $u \equiv 0$), we check for two things in evaluating the PASS method. Firstly, we see if PASS gives $u \equiv 1$ and $u \equiv 0$ for all possible parameters. Secondly, we check the number and value of false positives given by PASS. We present Tables 1(a) and 1(b) as results:

\begin{table}
\begin{center}

\begin{tabular}{ | c | c | }
\hline
$n$ & Number of Steady States \\
\bhline
4 & 2 \\
\hline
6 & 2 \\
\hline
8 & 2 \\ 
\hline
16 & 2 \\
\hline
32 & 2 \\
\hline
\end{tabular}
\hspace{.5in}
\begin{tabular}{ | c | c | }
\hline
$\mu$ & Number of Steady States \\
\bhline
.2 & 8 \\
\hline
.5 & 4 \\ 
\hline
2 & 2 \\
\hline
5 & 2 \\
\hline
\end{tabular}

\caption{Number of steady states in terms of $n$ and $\mu$.}
\label{Fisher Statistics}
\end{center}
\end{table}

The two steady states reported for all values of $n$ listed above were $u_i \equiv 0$ and $u_i \equiv 1$ as desired.

For $n = 16$, we then tested $\mu = .2$, $\mu = .5$, $\mu = 2$, and $\mu = 5$ ($\mu = 1$ was used in \cref{Fisher Statistics}(a)). The results are in \cref{Fisher Statistics}(b). Once again, $u_i \equiv 0$ and $u_i \equiv 1$ were found for all parameters. All of the additional steady states found for $\mu = .2$ and $\mu = .5$ were also constant-valued. For $\mu = .2$ the addition constant stady states were $u_i \equiv .05$, .10, .9, .95, 1.05, and 1.1. For $\mu = .5$ we had $u_i \equiv .05$ and $1.05$. Observe that all of these values are clustered around 0 and 1.

\subsection{Benjamin-Bona-Mahony}\label{sec:bbm}

The Benjamin-Bona-Mahony equation \citep{BBM} models long surface gravity waves---waves with large wavelength are caused by gravity. The equation is

\begin{equation}
\label{bbm}
    u_t + u_x + uu_x - u_{xxt} = 0
\end{equation}

This equation was completely solved in 1979 by Olver \citep{olver_1979}. The discrete form for solving steady states is:

\begin{equation}
    \label{bbm_discrete}
    \frac{u_{i+1} - u_i}{h}(1 + u_i) = 0
\end{equation}

From the equation, we can see that all constant-valued solutions are steady states of the equation. For $n = 16$, $h = .05$, and $u_i \in [0,2]$ with a resolution of $.05$, the pixel array associated with finding all numerical steady states to this system is:

\begin{figure}[h]
\centering
\includegraphics[width=6cm]{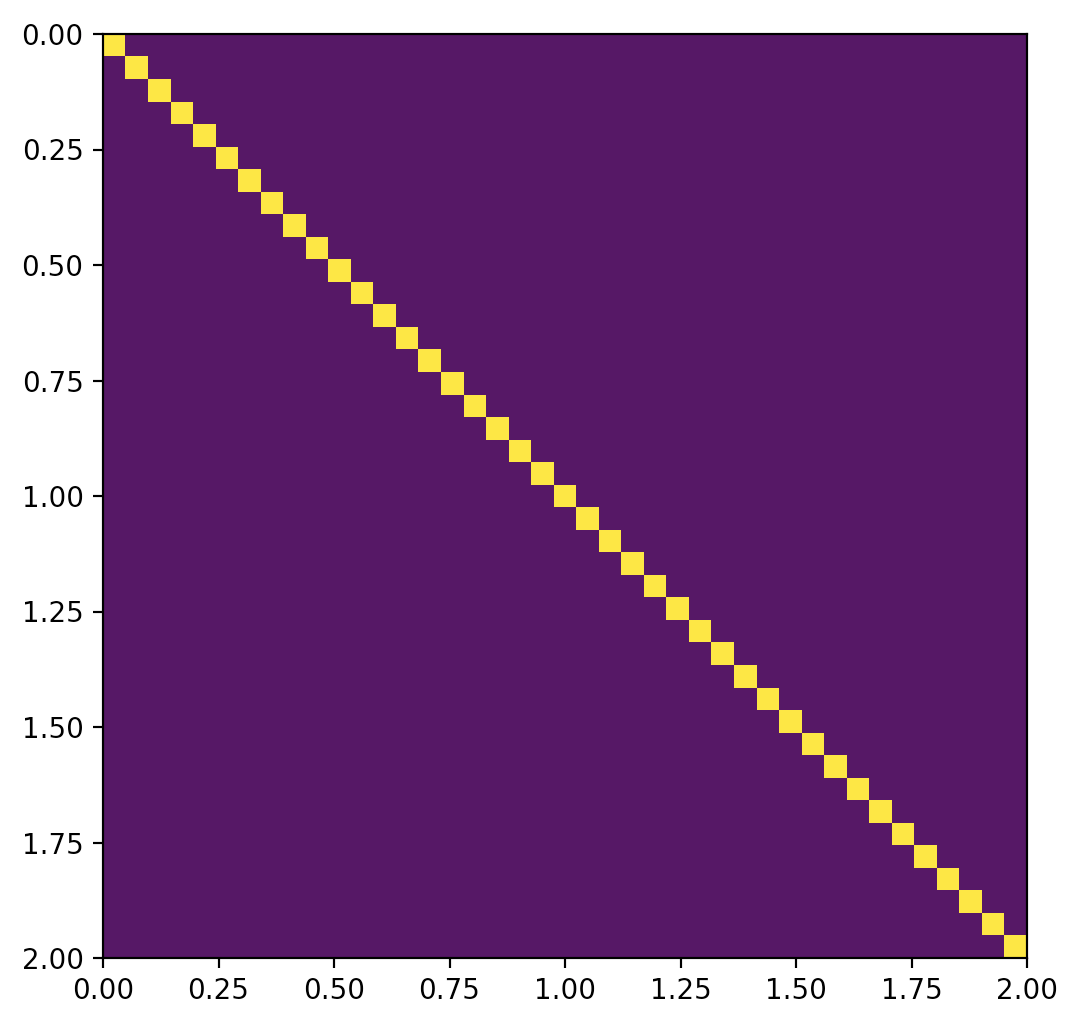}
\caption{The pixel array associated with the Benjamin-Bona-Mahony equation. Yellow indicates at least one steady state was found for the boundary conditions, purple otherwise}
\label{bbmPixel}
\end{figure}

One can verify from the output of the program that the numerical steady state produced for each boundary condition is indeed constant-valued and equal to the value at the boundaries.

\subsection{Sine-Gordon}\label{sec:sg}

Finally, we look at the Sine-Gordon equation, which is:

\begin{equation}
    \label{Sine}
    u_{tt} - u_{xx} + \sin u = 0
\end{equation}

Soliton solutions of this equation are $u = \pm 4 \tan^{-1} \exp(\pm \frac{x - Ut}{\sqrt{1 - U^2}})$ \citep{SG-behavior}. In the limit case of $t \rightarrow \pm \infty$, for $U \neq 0,1$, $u$ approaches the constant steady states $x = \pm 4 \pi/2 = \pm 2\pi$ or $x = 0$ (depending on the sign of $\frac{x-Ut}{\sqrt{1-U^2}}$). Both of these steady state solutions are numerical steady states of the discretized form of equation \ref{Sine}:

\begin{equation}
    \label{discrete SG}
    \sin u_i - \frac{u_{i+1} - 2u_i + u_{i-1}}{h^2} = 0
\end{equation}

However, discretization gives us false positive steady states. In fact, it gives us infinitely many false positives: $u_i \equiv k\pi$ for all integers $k$ $\forall i$. It is not obvious from observation whether there are even more numerical steady states. Thus, we run PASS on the discrete equation, and get, for $n = 8$, $b = .2$, $u_i \in [0,7]$:

\begin{figure}[h]
\centering
\includegraphics[width=6cm]{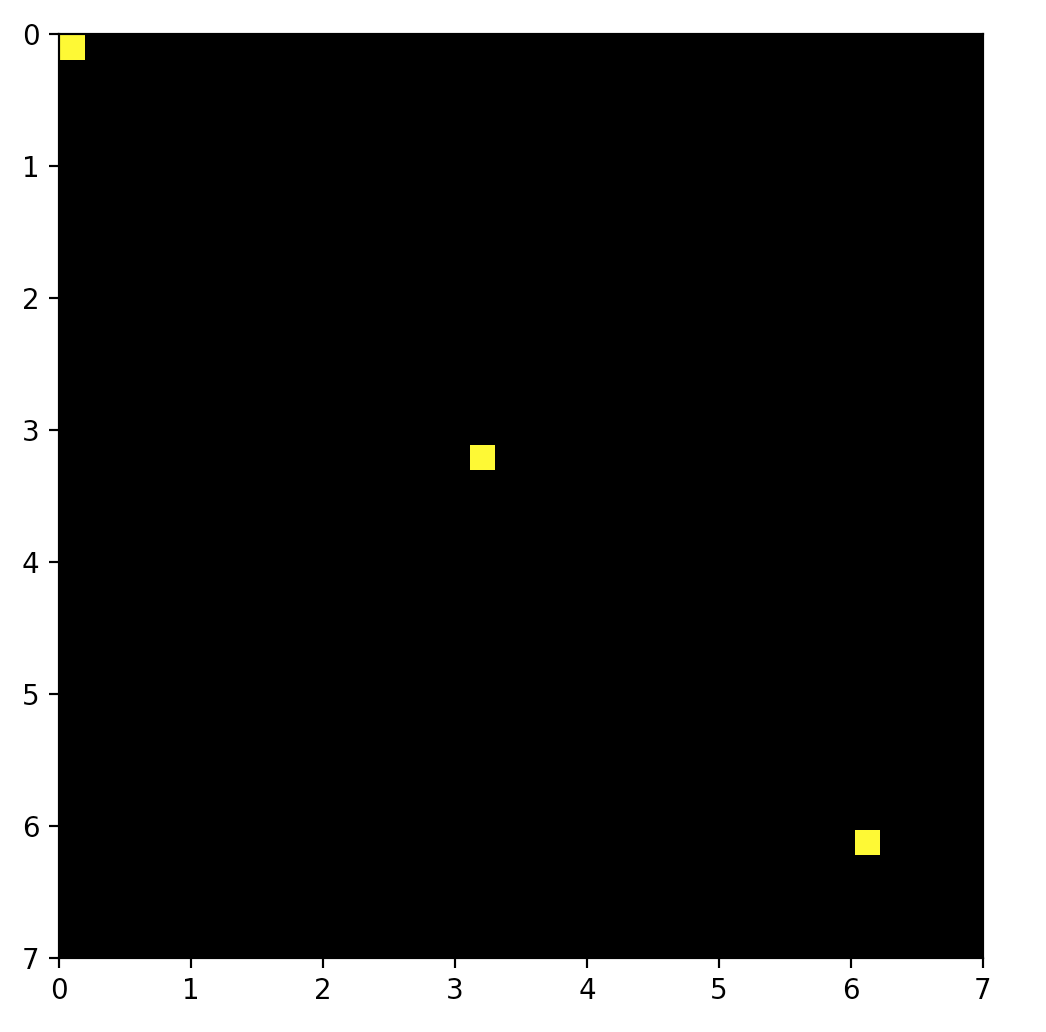}
\caption{The pixel array associated with Sine-Gordon equation.}
\label{sgPixel}
\end{figure}

From the pixel array, we see that there are numerical solutions indicated only for homogenous boundary values at 0 and $\approx \pi, 2\pi$. The output of the PASS program confirms that the steady states at those boundaries are constant.

\section{Discussion and Conclusion}

The PASS method successfully finds all steady states to various partial differential equations within a bounding box. If exact steady states cannot be found, approximate steady states that are similar to theoretical ones are found. Thus, PASS is effective at approximating numerical steady states of partial differential equations.

The primary novel benefit of the PASS method is the parallelizing the computation of many steady states, and the visualization of many steady states at once in the resultant pixel array. By doing so, it can facilitate the discovery and confirmation of patterns.

The PASS method has the same limitations as the original PA method, such as producing false positives and visualization of complex systems, and some additional limitations associated with discretization. As noted in section \cref{sec:sg}, discretization methods like finite difference may produce equations with more numerical steady states than the original. This issue may be resolved with more advanced discretization techniques such as finite element. 

Because PASS does a lot of array manipulation, it can be very slow, despite the original PA method being very efficient. This is because the product of solutionSets can result in an exponentially increasing number of steady states per pixel. Thus, even small systems analyzed with PASS might cause Python to crash if there are many steady states (such as with the heat equation). This problem can be mitigated by using PASS in conjunction with the PA method. In other words, because the PASS method is able to extract solution sets, it can be executed after running the more efficient boolean pixel array method to locate the areas of interest on which to run PASS.

The PASS method may be extended to find steady states of complex PDEs and even two- or three-dimensional systems. This may done by adjusting the steady state condition, the dimensions of our pixel arrays, editing plotting code for those arrays. These are changed based on the discretization and wiring diagram of whatever new equation we are examining. However, because many changes are required, using PASS might less practical for people compared to the usual PA method.

\section{Acknowledgments}

We thank our sponsors for the opportunity to do this research: CTL was sponsored by the MIT UROP program, and DIS was sponsored by AFOSR grants FA9550-14-1-0031 and FA9550-17-1-0058. We would also like to thank Professor Steven Johnson for suggesting this project to us.

\section{REFERENCES}

\bibliography{sources}

\begin{thebibliography}{8}
\providecommand{\natexlab}[1]{#1}
\providecommand{\url}[1]{\texttt{#1}}
\expandafter\ifx\csname urlstyle\endcsname\relax
  \providecommand{\doi}[1]{doi: #1}\else
  \providecommand{\doi}{doi: \begingroup \urlstyle{rm}\Url}\fi

\bibitem[BBM(2016)]{BBM}
Benjamin–bona–mahony equation, 2016.
\newblock URL
  \url{https://en.wikipedia.org/wiki/Benjamin-Bona-Mahony_equation}.

\bibitem[Baeumer et~al.(2008)Baeumer, Kovács, and Meerschaert]{Numerical_RD_1}
Boris Baeumer, Mihály Kovács, and Mark~M. Meerschaert.
\newblock Numerical solutions for fractional reaction–diffusion equations.
\newblock \emph{Computers \& Mathematics with Applications}, 55\penalty0
  (10):\penalty0 2212 -- 2226, 2008.
\newblock ISSN 0898-1221.
\newblock \doi{http://dx.doi.org/10.1016/j.camwa.2007.11.012}.
\newblock URL
  \url{http://www.sciencedirect.com/science/article/pii/S0898122107007420}.
\newblock Advanced Numerical Algorithms for Large-Scale Computations.

\bibitem[Liu et~al.(1996)Liu, Allen, Kojouharov, and Chen]{Numerical_RD_2}
B.~Liu, {M. B.} Allen, H.~Kojouharov, and B.~Chen.
\newblock \emph{Finite-element solution of reaction-diffusion equations with
  advection}, volume~1, pages 3--12.
\newblock Computational Mechanics Publ, 1996.

\bibitem[Reitz(1981)]{Numerical_RD_3}
Rolf~D. Reitz.
\newblock A study of numerical methods for reaction-diffusion equations.
\newblock \emph{SIAM Journal on Scientific and Statistical Computing},
  2\penalty0 (1):\penalty0 95--106, 1981.
\newblock \doi{10.1137/0902008}.
\newblock URL \url{https://doi.org/10.1137/0902008}.

\bibitem[{Spivak} et~al.(2016){Spivak}, {Dobson}, {Kumari}, and
  {Wu}]{Introduction_to_PA}
D.~I. {Spivak}, M.~R.~C. {Dobson}, S.~{Kumari}, and L.~{Wu}.
\newblock {Pixel Arrays: A fast and elementary method for solving nonlinear
  systems}.
\newblock \emph{ArXiv e-prints}, August 2016.

\bibitem[Berestycki et~al.(2009)Berestycki, Nadin, Perthame, and
  Ryzhik]{Fisher}
Henri Berestycki, Gregoire Nadin, Benoit Perthame, and Lenya Ryzhik.
\newblock The non-local fisher-kpp equation: Travelling waves and steady
  states.
\newblock 22:\penalty0 2813, 10 2009.

\bibitem[Olver(1979)]{olver_1979}
Peter~J. Olver.
\newblock Euler operators and conservation laws of the bbm equation.
\newblock \emph{Mathematical Proceedings of the Cambridge Philosophical
  Society}, 85\penalty0 (1):\penalty0 143–160, 1979.
\newblock \doi{10.1017/S0305004100055572}.

\bibitem[N.F.~Smyth(1999)]{SG-behavior}
A.L.~Worthy N.F.~Smyth.
\newblock Soliton evolution and radiation loss for the sine-gordon equation.
\newblock \emph{Physical Review}, 60\penalty0 (2):\penalty0 2330--2336, 1999.
\newblock URL
  \url{https://pdfs.semanticscholar.org/2c8a/d59e3bbec51ff2444da87313df817553cb2c.pdf}.

\end{thebibliography}

\end{document}